\documentclass[12pt]{article}
\usepackage{amssymb}
\begin{document}
\newtheorem{theorem}{Theorem}
\newtheorem{definition}{Definition}
\newtheorem{lemma}{Lemma}
\newtheorem{proposition}{Proposition}
\begin{center}{\Large $p$-Dirac Operators}\end{center}
\begin{center}{\Large Craig A. Nolder}\end{center}
\begin{center}{\Large Department of Mathematics, Florida State University, Tallahassee, Florida 32306-4510, USA}\end{center}
\begin{center}{and}\end{center}
\begin{center}{\Large John Ryan}\end{center}
\begin{center}{\Large Department of Mathematics, University of Arkansas, Fayetteville, AR 72701, USA}\end{center}
\begin{abstract} We introduce non-linear Dirac operators in $\mathbb{R}^{n}$ associated to the $p$-harmonic equation and we extend to other contexts including spin manifolds and the sphere.
\end{abstract}
\section{Introduction}
Associated to each type of Laplacian one usually sees a first order linearization, to a Dirac operator. For instance associated to the Laplacian in $\mathbb{R}^{n}$ is the euclidean Dirac operator arising in Clifford analysis. For the Laplace-Beltrami operator associated to a Riemannian manifold there is the Hodge-Dirac operator $d+ d^{\star}$, where $d$ is the exterior derivative and $d^{\star}$ is the Hodge codifferential which is the formal adjoint to $d$. Further, in reverse order, to the Atiyah-Singer-Dirac operator on a spin manifold there is the spinorial Laplacian. Also on $S^{n}$ one has a conformal Dirac operator and the conformal Laplacian. See for instance \cite{be,bds,lm,lr} for details.

\ Besides the Laplacian in $\mathbb{R}^{n}$ there are also the non-linear operators referred to as $p$-Laplacians. See for instance \cite{hkm,im,l}. Despite being non-linear these second order operators posses properties very similar to the usual Laplacian in euclidean space. Further when $p=2$ this operator corresponds to the usual Laplacian in euclidean space.

\ Here we shall introduce a first order nonlinear differential operator which in the case $p=2$ coincides with the euclidean Dirac operator. The conformal covariance of these operators are established. The $n$-harmonic equation arising here is Clifford algebra valued and the invariance of weak solutions to this equation under conformal transformations is only an invariance up to a factor of the pin group, the double covering of the orthogonal group. This is in contrast to weak solutions to the usual $n$-harmonic equation. We illustrate that we have a prpoer covariance not involving the pin group when we restrict to the scalar part of our clifford valued equations. 

\ Further a non-linear Cauchy-Riemann equation is introduced and its covariance under composition with non-constant holomorphic functions is described. 

\ Also a $p$-Dirac and a $p$-harmonic equation are set up on spin manifolds. We describe the behaviour of weak solutions to the $n$-Dirac equation under conformal rescaling of the metric on a spin manifold. 

\ We conclude by introducing $p$-Dirac and $p$-harmonic equations on the sphere $S^{n}$ and introducing solutions to these equations.
\newline
{\bf{Dedication}} This paper is dedicated to the memory of J. Bures.

\section{Preliminaries}
The $p$-Laplace equation is the non-linear differential equation $div\|\nabla f\|^{p-2}\nabla f=0$, where $f$ is a a sufficently smooth, scalar valued function defined on a domain in $\mathbb{R}^{n}$. Further the operator $\nabla$ is one of the simplest examples of a Dirac operator. One role here is to see how introducing a Dirac operator to the setting of $p$-Laplace equations might deepen ones perspective of such an equation. In order to introduce Dirac operators we need to first look at some basics of Clifford algebras.

\ Following \cite{p} and elsewhere one can consider $\mathbb{R}^{n}$ as embedded in the real Clifford algebra $Cl_{n}$. For each $x\in\mathbb{R}^{n}$ we have within $Cl_{n}$ the multiplication formula $x^{2}=-\|x\|^{2}$. If $e_{1},\ldots,e_{n}$ is an orthonormal basis for $\mathbb{R}^{n}$ this relationship defines an anti-commuting relationship $e_{i}e_{j}+e_{j}e_{i}=-2\delta_{ij}$. If this relationship is the only relationship assumed on $Cl_{n}$ then $1,e_{1},\ldots,e_{n},\ldots, e_{j_{1}}\ldots e_{j_{r}},\ldots,e_{1}\ldots e_{n}$ is a basis for $Cl_{n}$. Here $1\leq r\leq n$ and $j_{1}<\ldots j_{r}$. It follows that the dimension of $Cl_{n}$ is $2^{n}$. Further this algebra is associative.

\ We shall need the following antiautomorphism
\[\sim:Cl_{n}\rightarrow Cl_{n}:\sim(e_{j_{1}}\ldots e_{j_{r}})=e_{j_{r}}\ldots e_{j_{1}}.\]
For each $A\in Cl_{n}$ we shall write $\tilde{A}$ for $\sim(A)$.

\ Note for $x=x_{1}e_{1}+x_{2}e_{2}+\ldots +x_{n}e_{n}$ that $e_{1}(x)e_{1}=-x_{1}e_{1}+x_{2}e_{2}+\ldots +x_{n}e_{n}$. So we have a reflection in the $e_{1}$ direction. Similarly for $y\in S^{n-1}$, the unit sphere in $\mathbb{R}^{n}$, one has that $yxy$ is a reflection in the $y$ direction. Consequently for $y_{1},\ldots y_{J}\in S^{n-1}$ we have that $y_{1}\ldots y_{J}xy_{J}\ldots y_{1}$ is an orthogonal transformation acting on the vector $x$. In fact we have the group $Pin(n):=\{a\in Cl_{n}:a=y_{1}\ldots y_{J}:y_{1},\ldots y_{J}\in S^{n-1}$ and $J=1,2, 3,\ldots\}$. In \cite{p} and elsewhere it is shown that $Pin(n)$ is a double covering of the orthogonal group $O(n)$. When we restrict $J$ to be even we obtain a subgroup known as the spin group and denoted by $Spin(n)$. Further $Spin(n)$ is a double covering of the special orthogonal group, $SO(n)$. We shall also need the Lipschitz group $L(n)=\{a=x_{1}\ldots x_{J}:x_{1},\ldots, x_{J}\in \mathbb{R}^{n}\backslash\{0\}$ and $J\in\mathbb{N}\}$.

\ For $A=a_{0}+a_{1}e_{1}+\ldots+a_{1\ldots n}e_{1}\ldots e_{n}\in Cl_{n}$ we define the norm, $\|A\|$ , of $A$ to be $(a_{0}^{2}+\ldots+a_{1\ldots n}^{2})^{\frac{1}{2}}$. Conjugation on the Clifford algebra is defined to be the anti-automorphism $-:Cl_{n}\rightarrow Cl_{n}:-(e_{j_{1}}\ldots e_{j_{r}})=(-1)^{r}e_{j_{r}}\ldots e_{j_{1}}$. For $A\in Cl_{n}$ we write $\overline{A}$ for $-(A)$. Note that the real part, $Sc(A\overline{A})$, of $A\overline{A}$ is $\|A\|^{2}$. Further for $A$ and $B\in Cl_{n}$ the product $\overline{A}B$ defines a Clifford algebra valued inner product on $Cl_{n}$ for which $Sc(\overline{A}B)$ is the standard dot product on $\mathbb{R}^{2^{n}}$.

\ It is well known, see \cite{p}, that as a vector space $Cl_{n}$ is canonically isomorphic to the alternating al
gebra $\Lambda(\mathbb{R}^{n})$.

\ To the vector $x\in\mathbb{R}^{n}$ we can associate the differential operator $D:=\sum_{j=1}^{n}e_{j}\frac{\partial}{\partial x_{j}}$. This is the Dirac operator in euclidean space. Note that if $f$ is a $C^{1}$ real valued function defined on a domain $U$ in $\mathbb{R}^{n}$ then $Df=\nabla f$. Further $D^{2}=-\triangle_{n}$ where $\triangle_{n}$ is the Laplacian in $\mathbb{R}^{n}$.

\ In \cite{a} it is shown that given a M\"{o}bius transformation $M(x)$ over the one point compactification of $\mathbb{R}^{n}$ one can write this transformation as $(ax+b)(cx+d)^{-1}$ where $a$, $b$, $c$ and $d\in Cl_{n}$ and they satisfy the following conditions
\newline
(i) $a$, $b$, $c$ and $d$ are all products of vectors
\newline
(ii) $\tilde{a}c$, $\tilde{c}d$, $\tilde{d}b$ and $\tilde{b}a\in\mathbb{R}^{n}$
\newline
(iii) $\tilde{a}d-\tilde{b}c=1$.

\ If $y=M(x)$ then, \cite{r}, we have $cx+d\in L(n)$. Consequently $cx+d$ has a multiplicative inverse in $Cl_{n}$. It is shown in \cite{b} that $J_{-1}(M,x)^{-1}D_{x}J_{1}(M,x)=D_{y}$ where $D_{x}$ is the Dirac operator with respect to $x$ and $D_{y}$ is the Dirac operator with respect to $y$. Further $J_{-1}(M,x)=\frac{\widetilde{cx+d}}{\|cx+d\|^{n+2}}$ and $J_{1}(M,x)=\frac{\widetilde{cx+d}}{\|cx+d\|^{n}}$. Moreover $DJ_{1}(M,x)=0$. See \cite{r}. Consequently we have:

\begin{lemma}
Suppose $\psi$ is a $C^{1}$ function with compact support and $y=M(x)$. Then $D_{y}\psi(y)=(\widetilde{cx+d})^{-1}D_{x}(\widetilde{cx+d})\psi(M(x))$.
\end{lemma}
{\bf{Proof}} We know that $D_{y}\psi(y)=J_{-1}(M,x)^{-1}D_{x}J_{1}(M,x)\psi(M(x))$. But $DJ_{1}(M,x)=0$. The result now follows from Leibniz rule. $\Box$

\ It may be seen that $(\widetilde{cx+d})^{-1}D_{x}(\widetilde{cx+d})$ is a dilation and orthogonal transformation acting on $D_{x}$.

\section{$n$-Dirac and $n$-Laplace Equations and Conformal Symmetry}

\ If $v(x)$ is a $C^{1}$ vector field then the real or scalar part of $Dv$ is $div v(x)$. Keeping this in mind we formally define the $n$-Dirac equation for a $C^{1}$ function $f:U\rightarrow Cl_{n}$, with $U$ a domain in $\mathbb{R}^{n}$, to be $D\|f\|^{n-2}f=0$. This is a non-linear first order differential equation for $n>2$. When $f=Dg$ for some $Cl_{n}$ valued function $g$ then the $n$-Dirac equation becomes $D\|Dg\|^{n-2}Dg=0$. Further when $g$ is a real valued function the scalar part of this equation becomes $div(\|\nabla g\|^{n-2}\nabla g)=0$ which is the $n$-Laplace equation described earlier. In the Clifford algebra context the $n$-Laplace equation extends to the equation $D\|Du\|^{n-2}Du=0$. We shall refer to this equation as the $n$-$Cl_{n}$ Laplace equation. The function $ln\|x\|$ is a solution to this equation on $\mathbb{R}^{n}\backslash\{0\}$. When $u$ is scalar valued on identifying the Clifford algebra $Cl_{n}$ with the alternating algebra $\Lambda(\mathbb{R}^{n})$, the nonscalar part of the equation $D\|Du\|^{n-2}Du=0$ becomes $d\|du\|^{n-2}du=0$ where $d$ is the exterior derivative. When $n=3$ using the Hodge star map this equation becomes in vector calculus terminology $\nabla\times\|\nabla u\|^{n-2}\nabla u=0$.

\ Noting that $D\frac{x}{\|x\|^{n}}=0$ one may see that $\frac{x}{\|x\|^{2}}$ is a solution to the $n$-Dirac equation on $\mathbb{R}^{n}\backslash\{0\}$. We will assume that all $Cl_{n}$ valued test functions have components in $C^{\infty}_0(U)$.

\begin{definition}
Suppose $f:U\rightarrow Cl_{n}$ is in $L_{loc}^{n}(U)$, so each component of $f$ is in $L_{loc}^{n}(U)$. Then $f$ is said to be a weak solution to the $n$-Dirac equation if for each $Cl_{n}$ valued test function $\eta$ defined on $U$
\[\int_{U}(\overline{\|f\|^{n-2}f}D\eta)dx^{n}=0.\]
\end{definition}

\ Note that for $g\in W_{loc}^{1,n}(U)$ then $Dg$ is a weak solution to the $n$-Dirac equation.

\ We now proceed to establish a conformal covariance for the $n$-Dirac equation.

\begin{theorem}
Suppose $f:U\rightarrow Cl_{n}$ is a weak solution to the $n$-Dirac equation. Suppose also $y=M(x)=(ax+b)(cx+d)^{-1}$ is a M\"{o}bius transformation such that $cx+d$ is non-zero on the closure of $M^{-1}(U)$. Then $(cx+d)^{-1}f(M(x))$ is a weak solution to the $n$-Dirac equation on $M^{-1}(U)$.
\end{theorem}
{\bf{Proof}}: Consider $\int_{U}(\overline{\|f(y)\|^{n-2}f(y)}D_{y}\eta(y))dy^{n}$. As the Jacobian of $M$ is $\frac{1}{\|cx+d\|^{2n}}$ and $D_{y}=J_{-1}(M,x)^{-1}D_{x}J_{1}(M,x)$ this integral transforms to
\[\int_{M^{-1}(U)}(\overline{\|f(M(x))\|^{n-2}f(M(x))}J_{1}(M,x)D_{x}J_{1}(M,x)\eta(M(x))dx^{n}.\]
 Redistributing terms in $J_{1}(M,x)$ this integral becomes
\[\int_{M^{-1}(U)}(\overline{\|(cx+d)^{-1}f(M(x))\|^{n-2}(cx+d)^{-1}f(M(x))}D_{x}J_{1}(M,x)\eta(M(x)))dx^{n}.\]
As $cx+d$ is bounded on $M^{-1}(U)$ then $J_{1}(M,x)\eta(M(x))$ is a test function on $M^{-1}(U)$. Further as $cx+d$ is bounded and $C^{\infty}$ on $M^{-1}(U)$ then $(cx+d)^{-1}$ is a bounded $C^{\infty}$ function on $M^{-1}(U)$. Consequently $(cx+d)^{-1}f(M(x))\in L_{loc}^{n}(M^{-1}(U))$. The result follows. $\Box$

\begin{definition}
Suppose $f:U\rightarrow Cl_{n}$ belongs to $W_{loc}^{1,n}(U)$ and
\[\int_{U}(\overline{\|Df\|^{n-2}Df}D\eta)dx^{n}=0\]
for each $Cl_{n}$ valued test function defined on $U$. Then $f$ is called a weak solution to the $n$-$Cl_{n}$ Laplace equation.
\end{definition}

\ We shall now examine the conformal symmetry of weak solutions to the $n$-$Cl_{n}$ Laplace equation. Our arguments follow the lines for $A$-harmonic morphisms given in \cite{hkm}.

\begin{theorem}
Suppose $f:U\rightarrow Cl_{n}$ is a weak solution to the $n$-$Cl_{n}$ Laplace equation. Suppose further that $y=M(x)=(ax+b)(cx+d)^{-1}$ is a M\"{o}bius transformation and $cx+d$ is non-zero on the closure of $M^{-1}(U)$. Then $f(M(x))$ is a weak solution to the equation $D_{M}\|Df(M(x))\|^{n-2}D_{M}f(M(x))=0$ on $M^{-1}(U)$, where $D_{M}:=\Sigma_{j=1}^{n}\frac{(cx+d)}{\|cx+d\|}e_{j}\frac{\widetilde{cx+d}}{\|cx+d\|}\frac{\partial}{\partial x_{j}}$.
\end{theorem}
{\bf{Proof}} As $DJ_{1}(M,x)=0$ then on changing variables and applying Lemma 1 the integral $\int_{U}(\overline{\|Df\|^{n-2}Df}D\eta)dy^{n}$ becomes
\[\int_{M^{-1}(U)}\|cx+d\|^{2n}(\|D_{M}f(M(x))\|^{n-2}\overline{D_{M}f(M(x))}
D_{M}\eta(M(x)))\frac{dx^{n}}{\|cx+d\|^{2n}}\]
\[=\int_{M^{-1}(U)}\|D_{M}f(M(x))\|^{n-2}\overline{D_{M}f(M(x))}D_{M}\eta(M(x))dx^{n}.\]
In \cite{br} it is shown that $\|(cx+d)A\|=\|cx+d\|\|A\|$ for any $A\in Cl_{n}$. Consequently $\|D_{M}f(M(x))\|=\|Df(M(x))\|$. The result follows. $\Box$

\ Note that $\frac{\widetilde{cx+d}}{\|cx+d\|}$ belongs to the pin group $Pin(n)$. So the covariance we have described here for weak solutions to the $n$-harmonic equation is not the same as for the classical $n$-harmonic equation descibed in \cite{l} and elsewhere. We shall return to this point in the next section. First though let us note that it follows from \cite{br} that $Sc(D_{M}\|Df(M(x))\|^{n-2}D_{M}f(M(x)))=Sc(D\|Df(M(x))\|^{n-2}Df(M(x)))$ and 
\[Sc(\|Df(M(x)\|^{n-2}\overline{D_{M}f(M(x))}D_{M}\eta(M(x)))=Sc(\|Df(M(x))\|^{n-2}\overline{Df(M(x))}D\eta(M(x))).\]
 When $f$ is scalar valued this establishes the conformal invariance of the $n$-Laplace equation.

\section{$p$-Dirac and $p$-$Cl_{n}$ Laplace Equations and M\"{o}bius Transformations}

\ We now turn to the more general case. For any real positive number $p$ a differentiable function $f:U\rightarrow Cl_{n}$ is said to be a solution to the $p$-Dirac equation if $D\|f\|^{p-2}f=0$. For $1<p<n$ the function $\frac{x}{\|x\|^{\frac{n+p-2}{p-1}}}$ is a solution to this equation on $\mathbb{R}^{n}\backslash\{0\}$. We obtain this solution
by again noting that $D\frac{x}{\|x\|^{n}}=0$ and solving the equation $\|f\|^{p-2}f=\frac{x}{\|x\|^{n}}$.

\begin{definition}
Suppose that $f:U\rightarrow Cl_{n}$ belongs to $L_{loc}^{p}(U)$. Then $f$ is a weak solution to the $p$-Dirac equation if for each $Cl_{n}$ valued test function $\eta$ defined on $U$ we have $\int_{U}(\overline{\|f\|^{p-2}f}D\eta)dx^{n}=0$.
\end{definition}

\ Besides the $p$-Dirac equation we also need the following equation
\[D\|g\|^{p-2}A(x)g(x)=0\]
 where $g:U\rightarrow  Cl_{n}$ is a differentiable function and $A(x)$ is a real valued, smooth function.
We shall call this equation the $A,p$-Dirac equation. The $A,p$-Dirac equation is a natural generalization of the $A$-harmonic functions defined in \cite{hkm} and elsewhere.

\begin{definition}
Suppose that $g:U\rightarrow Cl_{n}$ is in $L_{loc}^{p}(U)$ and $A:U\rightarrow \mathbb{R}^{+}$ is a smooth bounded function. Then $g$ is a weak solution to the $A,p$-Dirac equation if for each $Cl_{n}$ valued test function $\eta$ defined on $U$ we have
\[\int_{U}(\overline{A(x)\|g(x)\|^{p-2}g(x)}D\eta(x))dx^{n}=0.\]
\end{definition}

\ By similar arguments to those used to prove Theorem 1 we now have:

\begin{theorem}
Suppose $g:U\rightarrow Cl_{n}$ is a weak solution of the $p$-Dirac equation and $y=M(x)=(ax+b)(cx+d)^{-1}$ is a M\"{o}bius transformation with $cx+d$ non-zero on the closure of $M^{-1}(U)$. Then $(cx+d)^{-1}g(M(x))$ is a weak solution to the $A,p$-Dirac equation on $M^{-1}(U)$, with $A(x)=\|cx+d\|^{p-n}$.
\end{theorem}

\begin{definition}
Suppose $h:U\rightarrow Cl_{n}$ is a solution to the equation
\[D\|Dh(x)\|^{p-2}Dh(x)=0\]
then $h$ is called a $p$-harmonic function.
\end{definition}

\ For $1<p<n$ the function $\|x\|^{\frac{p-n}{p-1}}$ is a solution to this equation. Again when $u$ is scalar valued one may identify $Cl_{n}$ with $\Lambda(\mathbb{R}^{n})$. In this case the non-scalar part of this $p$-harmonic equation becomes $d\|d u\|^{p-2}d u=0$. Also when $n=3$ the Hodge star map may be used to see that this equation becomes $\nabla\times\|\nabla u\|^{p-2}\nabla u=0$.

\ Note that when $h$ is real valued then the real part of the equation appearing in Definition 5 is the usual $p$-harmonic equation described in \cite{hkm}.

\begin{definition}
For a function $h:U\rightarrow Cl_{n}$ in $W_{loc}^{1,p}(U)$, then $h$ is called a weak solution to the $p$-harmonic equation if for each test function $\eta:U\rightarrow C   l_{n}$
\[\int_{U}(\overline{\|Dh\|^{p-2}Dh}D\eta)dx^{n}=0.\]
\end{definition}

\begin{definition}
For $h:U\rightarrow Cl_{n}$ a differentiable function and $A$ as in Definition 4, then $h$ is called an $A,p$-harmonic function if
\[DA(x)\|Dh(x)\|^{p-2}Dh(x)=0.\]
Further if $M(x)=(ax+b)(cx+d)^{-1}$ is a M\"{o}bius transformation then $h$ is called an $A,p,M$-harmonic function if 
\[D_{M}A(x)\|Dh(M(x))\|^{p-2}D_{M}h(M(x))=0.\]
\end{definition}

\begin{definition}
Suppose $h:U\rightarrow Cl_{n}$ belongs to $W_{loc}^{1,p}(U)$ and $A$ is as in Definition 4. Then $h$ is called a weak solution to the $A,p$-Laplace equation, or $A$, $p$-harmonic equation if for each test function $\eta:U\rightarrow Cl_{n}$
\[\int_{U}(\overline{A(x)\|Dh(x)\|^{p-2}Dh(x)}\eta(x))dx^{n}=0.\]
Further it is a weak solution to the $A,p,M$-harmonic equation if
\[\int_{M^{-1}(U)}\overline{A(x)\|Dh(M(x))\|^{p-2}D_{M}h(M(x))}D_{M}\eta(M(x))dx^{n}=0.\] 
\end{definition}

\ Further by similar arguments to those used to prove Theorem 2 we have:

\begin{theorem}
Suppose that $h:U\rightarrow Cl_{n}$ is a weak solution to the $p$-harmonic equation and $M(x)=(ax+b)(cx+d)^{-1}$ is a M\"{o}bius transformation with $cx+d$ non-zero on the closure of $M^{-1}(U)$. Then $h(M(x))$ is a weak solution to the $A,p,M$-harmonic equation on $M^{-1}(U)$ where $A(x)=\|cx+d\|^{2(p+2-n)}$.
\end{theorem}

Again $Sc(D_{M}A(x)\|Dh(M(x))\|^{p-2}D_{M}h(M(X)))=Sc(DA(x)\|Dh(M(x))\|^{p-2}Dh(M(x)))$ and 
\[Sc(A(x)\|Dh(M(x))\|^{p-2}\overline{D_{M}h(M(x))}D_{M}\eta(M(x)))=Sc(A(x)\|Dh(M(x))\|^{p-2}\overline{Dh(M(x))}D\eta(M(X))).\]
So when $h$ is scalar valued this again re-establishes the $A,p$ covariance of the $p$-harmonic equation.

\section{The $p$-Cauchy-Riemann Equation}

\ So far we have considered $p$-Dirac equations in dimensions $n\geq 3$. We now turn to look at the case $n=2$. In this setting the Dirac operator is $e_{1}\frac{\partial}{\partial x}+e_{2}\frac{\partial}{\partial y}$. This can be written as $e_{1}(\frac{\partial}{\partial x}+e_{1}^{-1}e_{2}\frac{\partial}{\partial y})$ and $\frac{\partial}{\partial x}+e_{1}^{-1}e_{2}\frac{\partial}{\partial y}=\frac{\partial}{\partial x}-e_{1}e_{2}\frac{\partial}{\partial y}=\frac{\partial}{\partial x}+e_{2}e_{1}\frac{\partial}{\partial y}$. Now $(e_{2}e_{1})^{2}=-1$. Consequently we can identify $e_{2}e_{1}$ with $i$, the square root of minus one. Then the operator $\frac{\partial}{\partial x}+e_{2}e_{1}\frac{\partial}{\partial y}$ can be identified with the Cauchy-Riemann operator $\frac{\partial}{\partial\overline{z}}$. If we restrict attention to functions taking values in the even subalgebra of $Cl_{2}$ spanned by $1$ and $e_{1}e_{2}$ and identify this algebra with $\mathbb{C}$ in the usual way then such a solution to the Dirac equation becomes a holomorphic function and vice versa.

\ A differentiable function $g:U\rightarrow\mathbb{C}$ is said to be a solution to the $p$-Cauchy-Riemann equation if it satisfies $\frac{\partial}{\partial\overline{z}}\|g(z)\|^{p-2}g(z)=0$. A function $g:U\rightarrow\mathbb{C}$ belonging to $L_{loc}^{p}(U)$ is said to be a weak solution to the $p$-Cauchy-Riemann equation if for each test function $\eta$ defined on $U$
\[\int_{U}\|g(z)\|^{p-2}\overline{g(z)}\frac{\partial}{\partial\overline{z}}\eta(z)dxdy=0.\]

\ Note that if $h:U\rightarrow\mathbb{C}$ is a $p$-harmonic function then $g(z):=\frac{\partial}{\partial z}h(z)$ is a solution to the $p$-Cauchy-Riemann equation.

\ Let us now suppose that $U$ is a bounded domain in the complex plane and $f:U\rightarrow\mathbb{C}$ is a non-constant holomorphic function with $f'(z)\ne 0$ on $U$. Using the identities $\eta(\zeta)=\frac{1}{\pi}\int_{U}\frac{\partial\eta(z)}{\partial\overline{z}}\frac{1}{z-\zeta}dxdy$ and $\eta(\zeta)=\frac{1}{\pi}\frac{\partial}{\partial\overline{z}}\int_{U}\frac{\eta(z)}{z-\zeta}dxdy$ for any test function $\eta:U\rightarrow\mathbb{C}$, and placing $z=f(u)$, then one may determine that
\[\frac{\partial}{\partial\overline{w}}f'(\zeta)^{-1}\eta(w)=\overline{f'}(\zeta)^{-1}\frac{\partial}{\partial\overline{\zeta}}\eta(f(\zeta))\]
where $w=f(\zeta)$.

\begin{theorem}
Suppose $g:U\rightarrow\mathbb{C}$ is a weak solution to the $p$-Cauchy-Riemann equation and $f(z)$ is a holomorphic function defined on $U$ with $f'(z)\ne 0$. Then $f'(\zeta)\|g(f(\zeta))\|^{p-2}g(f(\zeta))$ is a weak solution to the equation
\[\frac{\partial}{\partial\overline{\zeta}}f'(\zeta)\|g(f(\zeta))\|^{p-2}g(f(\zeta))=0.\]
\end{theorem}

\ The proof follows the same lines as the proof of Theorem 1.

\ Note that if $f'(\zeta)\|g(f(\zeta))\|^{p-2}g(f(\zeta))$ is differentiable, then as $f$ is holomorphic, $g(f(\zeta))$ is a solution to the $p$-Cauchy-Riemann equation.

\section{$p$-Dirac and $p$-Harmonic Sections on Spin Manifolds}

\ The material presented here  depends heavily on the dot product in $\mathbb{R}^{n}$. In fact one can readily extend many of the basic concepts given here to more general inner product spaces. We shall turn to the context of spin manifolds. Amongst other sources basic facts on spin manifolds can be found in \cite{lm}.

\ Suppose that $M$ is a connected, orientable, Riemannian manifold. Associated to such a manifold is a principle bundle with each fiber isomorphic to the group $SO(n)$. If this bundle has a lifting to a further principle bundle with each fiber isomorphic to the group $Spin(n)$, then $M$ is said to have a spin structure and $M$ is called a spin manifold. Associated to a spin manifold is a vector bundle $Cl(M)$ with each fiber isomorphic to $Cl_{n}$.

\ The Levi-Civita connection $\nabla$ on $M$ lifts to a connection $\nabla'$ on the spin structure. Associated to that connection is the Atiyah-Singer-Dirac operator $D'$. If $e_{1}(x),\ldots,e_{n}(x)$ is a local orthonormal basis on $M$ then locally $D'=\Sigma_{j=1}^{n}e_{j}(x)\nabla_{e_{j}(x)}$. Further \cite{er} the inner product associated to the Riemannian structure of $M$ lifts to a Clifford algebra valued inner product on $Cl(M)$. We denote this inner product by $<$ , $>$.

\ Suppose now that $U$ is a domain in $M$ and $f:U\rightarrow Cl(M)$ is a differentiable section. Then $f$ is said to be a solution to the $p$-Atiyah-Singer-Dirac equation if $D'\|f(x)\|^{p-2}f(x)=0$. Further a section $f:U\rightarrow Cl(M)$ belonging to $L_{loc}^{p}(U)$ is said to be a weak solution to the $p$-Atiyah-Singer-Dirac equation if for each test section $\eta:U\rightarrow Cl(M)$ we have $\int_{U}<\|f\|^{p-2}f,D'\eta>dU=0$ where $dU$ is the volume element induced by the metric on $M$.

\ Besides the $p$-Atiyah-Singer-Dirac equation we may also introduce $p$-spinorial harmonic functions. A twice differentiable section $h:U\rightarrow Cl(M)$ is said to be $p$-spinorial harmonic if $D'\|D'h\|^{p-2}D'h=0$. Further if we assumed that $h\in W_{loc}^{1,p}(U)$, then $h$ is a weak solution of the $p$-spinorial harmonic equation if for each test section $\eta:U\rightarrow Cl(M)$ we have
\[\int_{U}<\|D'h\|^{p-2}D'h,D'\eta>dU=0.\]

\ As $Cl_{n}$ contains an identity there is a projection operator $Sc:Cl(M)\rightarrow Cl_{\mathbb{R}}(M)$ where $Cl_{\mathbb{R}}(M)$ is the line bundle of $Cl(M)$ with each fiber the real part of the fiber of $Cl(M)$. It makes sense to now talk of the equation 
\begin{equation}
\int_{U}Sc<\|D'h\|^{p-2}D'h,D'\eta>dU=0.
\end{equation}

\ This last integral arises from the vanishing of the first variation associated to the Dirichlet integral $\int_{U}\|D'h\|^{p}dU$. When $p=2$ this integral gives rise to the spinorial Laplace equation $D'^{2}h=0$.

\ We can go a little further if our manifold is both a spin manifold and a conformally flat manifold. A manifold is said to be conformally flat if it has an atlas whose transition functions are M\"{o}bius transformations.

\ Suppose $M$ is a conformally flat spin manifold. For  a M\"{o}bius transition function $M:U\subset\mathbb{R}^{n}\rightarrow V\subset\mathbb{R}^{n}$ with $M(x)=(ax+b)(cx+d)^{-1}=(-ax-b)(-cx-d)^{-1}$ we can make an identification $(x,X)\leftrightarrow (M(x),\pm (cx+d)^{-1}X)$ where $x\in U$ and $X\in Cl_{n}$. As $M$ is a spin manifold  signs can be chosen so that these identifications are globally compatible over the manifold $M$. Consequently we have a vector bundle on $M$. Given the conformal covariance of the $n$-Dirac equation described in Theorem 1 it now follows from Theorem 1 that one can set up weak solutions to the $n$-Dirac equation over domains in $M$, and taking values in this vector bundle.

\ Similarly one can now use Theorem 2 and the remarks following it to see the conformal invariance of Equation 1.

\ Two metrics $g_{ij}$ and $g'_{ij}$ on a Riemannian manifold are said to be conformally equivalent if there is a function $k:M\rightarrow \mathbb{R}^{+}$ such that $g'_{ij}(x)=k(x)g_{ij}(x)$ for each $x\in M$. We now investigate how weak solutions to the $n$-Dirac equation transform under such conformal changes of metric on a spin manifold. We shall denote the inner product on the spinor bundle of $M$ associted to the metric $g_{ij}$ by $<$ , $>_{1}$ and the inner product associated to $g'_{ij}$ by $<$ , $>_{2}$. Further we denote the respective norms by $\|$ $\|_{1}$ and $\|$ $\|_{2}$. We denote the Dirac operator associated to $<$, $>_{1}$ by $D_{1}$ and the Dirac operator associated to $<$ , $>_{2}$ by $D_{2}$. Consequently the integral
\[\int_{U}<\|f\|_{2}^{n-2}f,D_{2}\eta>_{2}dU\]
becomes
\[\int_{U}<<f,f>_{2}^{\frac{n-2}{2}}f,D_{2}\eta>_{2}dU\]
\[=\int_{U}<<f(x),f(x)>_{2}^{\frac{n-2}{2}}f(x),D_{2}\eta(x)>_{1}k(x)^{2n}dU.\]
However, $D_{1}k^{n-1}(x)=k^{n+1}(x)D_{2}$. See for instance \cite{er}.

\ Consequently the previous integral becomes
\[\int_{U}<<f(x),f(x)>_{1}^{\frac{n-2}{2}}k^{n-2}(x)f(x),k(x)^{-n-1}(x)D_{1}k^{n-1}(x)\eta(x)>_{1}k^{2n}(x)dU.\]
This is equal to

\[\int_{U}<\|k(x)f(x)\|_{1}^{n-2}k(x)f(x),k(x)^{n-2}D_{1}k^{n-1}(x)\eta(x)>_{1}dU.\]

\ This calculation describes the change in the $n$-Dirac equation under conformal rescaling of the metric on a spin manifold. Similar transformations are possible for weak solutions to the $p$-Dirac equation under conformal changes in metric.

\section{The $p$-Dirac and $p$-Harmonic Equation on $S^{n}$}

\ Here we shall consider the unit sphere $S^{n}$ in $\mathbb{R}^{n+1}=span\{e_{1},\ldots,e_{n+1}\}$, and we shall consider functions defined on domains on $S^{n}$ and taking values in the Clifford algebra $Cl_{n+1}$. The stereographic projection from $S^{n}\backslash\{e_{n+1}\}$ to $\mathbb{R}^{n}$ corresponds to the Cayley transformation. Consequently one might expect that $p$-Dirac and $p$-harmonic equations can be set up on $S^{n}$. This indeed is the case. In \cite{lr} and elsewhere it is shown that the Dirac operator on $\mathbb{R}^{n}$ conformally transforms to the conformal Dirac operator $D_{S}:=x(\Gamma+\frac{n}{2})$ on $S^{n}$ where $\Gamma=\Sigma_{1\leq j<k\leq n}e_{i}e_{j}(x_{i}\frac{\partial}{\partial x_{j}}-x_{j}\frac{\partial}{\partial x_{i}})$ and $x\in S^{n}$.

\begin{definition}
Suppose $U$ is a domain on $S^{n}$ and $f:U\rightarrow Cl_{n+1}$ is a differentiable function. Then $f$ is called a solution to the $p$-spherical Dirac equation if $D_{S}\|f\|^{p-2}f=0$.
\end{definition}

\ Note that for $y\in S^{n}$ the function $\frac{x-y}{\|x-y\|^{2}}$ is a solution to the $n$-spherical Dirac equation. This follows as under the Cayley transformation the Clifford-Cauchy kernel $\frac{u-v}{\|u-v\|^{n}}$ in $\mathbb{R}^{n}$ conformally transforms to $\frac{x-y}{\|x-y\|^{n}}$ on $S^{n}$ and $D_{S}\frac{x-y}{\|x-y\|^{n}}=0$. For the same reason for $1<p\leq n$ the function $\frac{x-y}{\|x-y\|^{\frac{n+p-2}{p-1}}}$ is a solution to the $p$-spherical Dirac equation.

\begin{definition}
Suppose $U$ is a domain on $S^{n}$ and $f:U\rightarrow Cl_{n+1}$ belongs to $L^{p}(U)$. Then $f$ is a weak solution to the $p$-spherical Dirac equation if for each test function $\eta:U\rightarrow Cl_{n+1}$
\[\int_{U}(\overline{\|f\|^{p-2}f}D_S\eta)dU=0\]
where $dU$ is a volume element arising from the Lebesgue measure on $S^{n}$.
\end{definition}

\ One needs to be a bit careful in setting up a $p$-harmonic equation on the sphere. This is because the differential operator on $S^{n}$ that is conformally equivalent to the Laplacian in $\mathbb{R}^{n}$ is not $D_{S}^{2}$ but is the conformal Laplacian or Yamabe operator $Y_{S}$ described in \cite{be} and elsewhere. In \cite{br} it is shown that $Y_{S}=D_{S}(D_{S}-x)$.

\ In \cite{lr} it is shown that $(D_{S}+\frac{p}{2}x)\|x-y\|^{-n+p}=\frac{-n+p}{2}\frac{x-y}{\|x-y\|^{n-p}}$. Bearing this in mind and that the fundamental solution to $D_{S}$ is $\frac{x-y}{\|x-y\|^{n}}$ we define the $p$-spherical harmonic equation as follows.

\begin{definition}
Suppose $U$ is a domain on $S^{n}$ and $f:U\rightarrow Cl_{n}$ belongs to $W_{loc}^{1,p}(U)$. Then $f$ is a weak solution to the $p$-spherical harmonic equation if weakly $D_{S}\|(D_{S}+\frac{p}{2}x)f(x)\|^{p-2}(D_{S}+\frac{p}{2}x)f(x)=0$.
\end{definition}

\ Solutions to the $p$-spherical harmonic equation include $\|x-y\|^{\frac{p-n}{p-1}}$.

\end{document}